\documentclass[11pt]{amsart}
\usepackage{amssymb}
\usepackage{amsfonts}
\usepackage{color}
\usepackage{setspace}
\input{diagrams}

\newarrow{Dash}{}{dash}{}{dash}>
\newarrow{Into}C--->
\newarrow{Onto}----{>>}
\newarrow{Equal}===={=}
\newarrow{Dots}....>

\newtheorem{theorem}{Theorem}[section]
\newtheorem{proposition}[theorem]{Proposition}
\newtheorem{lemma}[theorem]{Lemma}
\newtheorem{corollary}[theorem]{Corollary}
\newtheorem{remark}[theorem]{Remark}

\renewenvironment{proof}{{\noindent\bf Proof.}}{\hfill $\Box$\par\vskip3mm}

\newcommand{\Cc}{\mathcal{C}}
\newcommand{\Dd}{\mathcal{D}}

\newcommand{\Ff}{\mathcal{F}}

\newcommand{\Ll}{\mathcal{L}}
\newcommand{\Mm}{\mathcal{M}}
\newcommand{\Nn}{\mathcal{N}}

\newcommand{\Rr}{\mathcal{R}}

\begin{document}
\title{Co-Frobenius Coalgebras}

\begin{abstract}
We investigate left and right co-Frobenius coalgebras and give equivalent characterizations which prove statements dual to the characterizations of Frobenius algebras. We prove that a coalgebra is left and right co-Frobenius if and only if $C\cong Rat(C^*_{C^*})$ as right $C^*$-modules and also that this is eqhivalent to the fact that the functors ${\rm Hom}_K(-,K)$ and ${\rm Hom}_{C^*}(-,C^*)$ from $\Mm^C$ to ${}_{C^*}M$ are isomorphic. This allows a definition of a left-right symmetric concept of co-Frobenius coalgebras that is perfectly dual to the one of Frobenius algebras and coincides to the existing notion left and right co-Frobenius coalgebra.
\end{abstract}

\author{Miodrag Cristian Iovanov}
\thanks{$^*$ Partially supported by the Flemish-Romanian project "New Techniques in Hopf algebras and Graded Rings Theory" (2005) and from the BD type grant BD86 2003-2005 by CNCSIS}
\date{}
\maketitle

\section*{Introduction}

Let $A$ be an algebra. $A$ is called Frobenius algebra if the right regular module $A$ is isomorphic to the right $A$ module $A^*$. This is known to be equivalent to the fact that the functors ${\rm Hom}_A(-,A)$ and ${\rm Hom}_K(-,K)$ from $\Mm_A$ to ${}_{A}\Mm$ are naturally isomorphic (see \cite{CR}). A functor $F:\Cc\rightarrow \Dd$ is called Frobenius if it has the same left and right adjoint (see \cite{CMZ}, \cite{NT1}). By using this concept, a Frobenius algebra can equivalently characterized by the fact that the forgetful functor $F$ from $\Mm_A$ to $M_K$ is a Frobenius functor. Moreover, an algebra $A$ is Frobenius if and only if there is a $K$ bilinear form $(-,-)$ on $A\times A$ that is associative and (left) non-degenerate. These definitions turn out to be left-right symmetric and imply the finite dimensionality of the algebra $A$. \\
Let $C$ be a coalgebra over a field $K$ and let $C^*$ be the dual coalgebra. The coalgebra $C$ is called {\emph right co-Frobenius} if there is a monomorphism of right $C^*$ modules from $C$ into $C^*$, or equivalently, there is a $C^*$ balanced bilinear form on $C\times C$ which is right non-degenerate. The concept of {\emph left co-Frobenius} coalgebra is defined by symmetry. A coalgebra is called simply co-Frobenius if it is both left and right co-Frobenius. The notion of right co-Frobenius coalgebra was introduced by Lin in \cite{L} and it abstracts a relevant coalgebra property of Hopf algebras with non-zero integral (\cite{L}, Theorem 3). A study of one side co-Frobenius coalgebras and more generally of one side quasi co-Frobenius coalgebras was carried out in \cite{NT1} and \cite{NT2}. For a Hopf algebra $H$, it is known that $H$ is right co-Frobenius is equivalent to $H$ having non-zero right integral, so this property defines this class of Hopf algebras. Hopf algebras with non-zero integrals are very important and have been intensely studied in the last years, mainly because they have very good structural and representation theoretic properties (see \cite{AD}, \cite{DNR}, \cite{H1}, \cite{S1}, \cite{S2}, \cite{Sw1}, \cite{Sw2} and references therein). Quantum groups with non-zero integrals are also of great interest; see \cite{APW}, \cite{AD}, \cite{H1}, \cite{H2}, \cite{H3}.\\
Lin proves that for a Hopf algebra $H$ being right co-Frobenius amounts to being left co-Frobenius and therefore it defines a left-right symmetric concept for Hopf algebras. However, that is not the case for coalgebras, as an example in \cite{L} shows that right co-Frobenius coalgebras need not be left co-Frobenius (see also \cite{DNR}, example 3.3.7). Therefore the question of whether a left-right symmetric concept of co-Frobenius coalgebra can be defined, a concept that would be dual to that of Frobenius algebras, would recover the notion of Frobenius algebra in the finite dimensional case and the notion co-Frobenius coalgebra. It is proved in \cite{NT2} that that a coalgebra $C$ is co-Frobenius if and only if $Rat({}_{C^*}C^*)\cong {}_{C^*}C$ and $Rat(C^*_{C^*})\cong C_{C^*}$. It is then natural to ask which coalgebras satisfy the property $Rat({}_{C^*}C^*)\cong {}_{C^*}C$ and if this property $Rat({}_{C^*}C^*)\cong {}_{C^*}C^*$ and if these are exactly the co-Frobenius coalgebras. In order to be able to make an analogue to the algebra case, we will consider the duality functors $C^*$-dual ${\rm Hom}_{C^*}(-,C^*)$ and $K$-dual ${\rm Hom}(-,K)$ from ${\mathcal M}^C$ to ${}_{C^*}{\mathcal M}$. For the dual of the right regular comodule $C$ it is natural to look at $C^*$ and take its rational part $Rat(C^*_{C^*})$. \\
The main result of the paper states that if $C^C$ is isomorphic to the dual $Rat({}_{C^*}C^*)$ of the left comodule ${}^{C}C$, then the coalgebra $C$ is co-Frobenius and thus this property defines the notion of co-Frobenius coalgebra. As a consequence, we also show that $C$ is co-Frobenius if and only if the two duality functors ${\rm Hom}_{C^*}(-,C^*)$ and ${\rm Hom}(-,K)$ are isomorphic when evaluated on right comodules and furthermore, this is also equivalent to the existence of a bilinear form on $C$ which is $C^*$-balanced and is left and right non-degenerate. These results generalize and extend known results from the algebra case; however, the proofs are completely different from the ones in the algebra case and involve the use of several techniques and results specific to coalgebra theory. 


\section{Preliminary results}
We recall some results on coalgebras, which we state with references for completeness of the text. For basic facts on coalgebra and comodule theory one should see \cite{A}, \cite{DNR}, \cite{M} and \cite{Sw1}. Let $C$ be a coalgebra. If $M$ is a finite dimensional right $C$-comodule, then it becomes left $C^*$-module by $c^*\cdot m=m_0c^*(m_1)$ and its dual $M^*={\rm Hom}(M,K)$ becomes a right $C^*$-module (as stated above) which is rational, so it is a left $C$-comodule. The following results are Proposition 4, page 34 from \cite{D} and Lemma 15 from \cite{L} respectively. See also Corollaries 2.4.19 and 2.4.20 from \cite{DNR}.

\begin{lemma}\label{i1}
Let $Q$ be a finite dimensional right $C$-comodule. Then $Q$ is injective (projective) as a left $C^*$-module if and only if it is injective (projective) as a right $C$-comodule.
\end{lemma}

\begin{lemma}\label{i2}
Let $M$ be a finite dimensional right $C$-comodule. Then $M$ is an injective right $C$-comodule if and only if $M^*$ is a projective left $C$-comodule.
\end{lemma}


Let $C$ be a coalgebra over a field $K$. Denote by $\Ll$ and respectively $\Rr$ a set of representatives for the types of isomorphism of simple left, respectively right $C$-comodules. Let $s({}^{C}C)=\bigoplus\limits_{i\in I}S_i$ be a decomposition of the left socle of $C$ into simple left comodules and $s(C^C)=\bigoplus\limits_{j\in J}T_j$ a decomposition of the right socle. Then we have a direct sum decomposition $C=\bigoplus\limits_{i\in I}E(S_i)$ as left $C$-comodules and $C=\bigoplus\limits_{j\in J}E(T_j)$ as right $C$-comodules, where $E(S_i)$ are injective envelopes of $S_i$ contained in $C$ and $E(T_j)$ are right comodules that are injective envelopes of $T_j$ contained in $C$. Also every simple left $C$-comodule is isomorphic to one of the $S_i$'s and similarly every right simple $C$-comodule is isomorphic to one of the $T_j$'s. Note that the sets $\Ll$, $\Rr$, $I$ and $J$ are of the same cardinality. See \cite{G} and \cite{DNR}.

Let $C$ be a coalgebra. Recall from \cite{L} that $C$ is said to be right co-Frobenius if there is a monomorphism of right $C^*$-modules from $C$ into $C^*$. The notion of left co-Frobenius alegebra is defined similarly. $C$ is called right (left) semiperfect if every finite dimensional right (left) comodule has a projective cover. Also recall from \cite{L} that $C$ is right semiperfect if and only if the injective envelope of any simple left comodule is finite dimensional, and that a right co-Frobenius coalgebra is also right semiperfect. \\
The following proposition shows that $C\cong Rat({}_{C^*}C^*)$ as left $C^*$-modules for a left and right co-Frobenius coalgebra. See Corollary 1.2 from \cite{CDN} for the proof. 


\begin{proposition}\label{i3}
If $C$ is a left and right co-Frobenius coalgebra, then any injective morphism of right $C$ modules $\varphi:C\rightarrow(Rat{}_{C^*}C^*)$ is an isomorphism.
\end{proposition}

\begin{lemma}\label{i4}
Let $S$ be a simple left comodule and $E(S)$ be an injective envelope of $S$ contained in $C$. Then $S^\perp=\{\alpha\in E(S)^*\mid \alpha\mid_{S}=0\}$ is a maximal and small left $C^*$ submodule of $E(S)^*$ and $E(S)^*$ is generated by any $f\notin S^\perp$. Consequently, $E(S)^*$ is an indecomposable left $C^*$-module.
\end{lemma}
\begin{proof}
We begin by the following remark: for any left comodule $M$, we have a left $C^*$-modules isomorphism ${\rm Hom}^C(M,C)\cong {\rm Hom}(M,K)$ given by $f\mapsto \varepsilon\circ f$, where $\varepsilon$ is the counit of $C$. Here the left $C^*$-module structure on ${\rm Hom}^C(M,C)$ is given by the left $C^*$-action on $C$, namely $(c^*\cdot f)(x)=c^*\cdot f(x)=c^*(f(x)_2)f(x)_1$. For $M=C$ this isomorphism becomes even an isomorphism of algebras $({\rm Hom}^C(C,C),+,\circ)\cong (C^*,+,*)$, where $\circ$ is the composition of morphisms and $*$ is the usual convolution on $C^*$. It is not difficult to see that by this isomorphism ${\rm Hom}^C(M,C)$ becomes a left ${\rm Hom}^C(C,C)$ module and that the structure is given simply by composition, namely for $\alpha\in{\rm Hom}^C(C,C)$ and $f\in{\rm Hom}^C(M,C)$, $(\alpha\rightarrow f)=\alpha\circ f$. Thus we may proof the statement equivalently for the left ${\rm Hom}^C(C,C)$ module ${\rm Hom}^C(M,C)$.\\
Let $S$ be a simple subcomodule of $C$ and $E(S)$ an injective envelope of $S$ contained in $C$. Then there is a left subcomodule $X$ of $C$ such that $E(S)\oplus X=C$ in ${}^{C}{\Mm}$. As the functor ${\rm Hom}^C(-,C)$ is exact, we obtain an epimorphism $\pi:{\rm Hom}(E(S),C)\longrightarrow {\rm Hom}(S,C)$, $\pi(f)=f\mid{}_{C^*}$. The kernel of this morphism coincides with $S^\perp$ 
Let $f\in {\rm Hom}^C(E(S),C)$ such that $f\mid_{S}\neq 0$. Then $Ker(f)\cap S=0$ as $S$ is simple and so $Ker\,f=0$ because $S$ is essential in $E(S)$. So $E(S)\cong f(E(S))$ and thus there is a left subcomodule $M$ of $C$ so that $C=f(E(C))\oplus M$. We can extend $f$ to a left comodule isomorphism $\overline{f}$ from $C$ to $C$ by taking $\overline{f}$ to be $f$ on $E(S)$ and denote by $h$ its inverse. Now if $g$ is another element of ${\rm Hom}^C(E(S),C)$, denoting by $\overline{g}$ its extension to $C$ that equals $0$ on $X$ and $g$ on $E(S)$ we have $\overline{g}=\overline{g}\circ id_{C}=\overline{g}\circ(h\circ\overline{f})=(\overline{g} \circ h)\circ\overline{f}$ and then restricting to $E(S)$, $g=(\overline{g}\circ h)\circ f$, equivalently $g=(\overline{g}\circ h)\rightarrow f$ (by the left action of ${\rm End}({}^{C}C)$ on ${\rm Hom}^C(E(S),C)$) showing that ${\rm Hom}^C(E(S),C)$ is generated by $f$, for an arbitrary $f\in {\rm Hom}^C(E(S),C)\setminus S^\perp$. This shows that if $M\subsetneq E(S)^*$ is a submodule of $E(S)^*$, then $M\subseteq S^\perp$, so $S^\perp$ is the only maximal subcomodule of $E(S)^*$, and also $S^\perp\ll E(S)^*$. Consequently if $E(S)^*=M\oplus N$, then if $M,N\neq E(S)^*$ we get $M,N\subset S^\perp$ which is a contradiction as $S^\perp\neq E(S)^*$ as $S\neq 0$.
\end{proof}

\section{The Main result}

\begin{proposition}\label{s1}
Suppose $C\cong Rat(C^*_{C^*})$. Then for every $T\in\Rr$, we have either $E(T)$ finite dimensional and $E(T)^*\cong E(S)$ for some $S\in \Ll$, or $E(T)$ is infinite dimensional and $Rat(E(T)^*)=0$. 
\end{proposition}
\begin{proof}
From the hypothesis, $C$ is right co-Frobenius and hence right semiperfect and then the $E(S_i)$'s are finite dimensional. Then we have $$\bigoplus\limits_{i\in I}E(S_i)\cong Rat(C^*_{C^*})=Rat(\prod\limits_{j\in J}E(T_j)^*)=Rat(E(T_k)^*)\oplus Rat(\prod\limits_{j\in J\setminus \{k\}}E(T_j)^*)$$ for any $k\in J$. Then $Rat(E(T_k)^*)$ is injective and isomorphic to a direct sum of $E(S_i)$. By Lemma \ref{i1} we have that the $E(S_i)$'s are injective in $\Mm_{C^*}$ and by Lemma \ref{i4} $E(T_k)^*$ is indecomposable. Then we have two possibilities:\\
-$Rat(E(T_k)^*)=0$ which implies that $E(T)$ is infinite dimensional (as otherwise $E(T)^*$ is finite dimensional and rational)\\
-$Rat(E(T_k))^*\neq 0$ and then there is a direct sum decomposition $\bigoplus\limits_{i\in I'}E(S_i)\cong Rat(E(T_k)^*)$ so there is an $i\in I$ such that $E(S_i)$ is a direct summand of $E(T_k)^*$ (because $E(S_i)$ is injective in $\Mm_{C^*}$) and then $E(S_i)\cong E(T_k)^*$ (because $E(T_k)^*$ is an indecomposable right $C^*$-module). As every $T\in \Rr$ is isomorphic to one of the $T_j$'s, the proposition is proved.
\end{proof}

Denote by $J_0=\{j\in J\mid Rat(E(T_j)^*)\neq 0\}$ and $J'=J\setminus J_0$. Notice that $Rat(\prod\limits_{j\in J'}E(T_j)^*)=0$. Indeed, denoting by $p_j$ the canonical projection on the $j$'th component of the direct product, we have that if $0\neq x\in Rat(\prod\limits_{j\in J'}E(T_j)^*)$, then there is $j\in J$ such that $p_j(x)\neq 0$. But $p_j(x)\in Rat(E(T_j)^*)=0$, which is a contradiction. 

\begin{corollary}\label{s2}
With the above notations we have ${}^{C}C\cong Rat(\prod\limits_{j\in J_0}E(T_j)^*)$ provided that $C\cong Rat(C^*_{C^*})$.
\end{corollary}

\begin{proposition}\label{s3}
If $C$ is right semiperfect then the set $\{E(S)^*\mid S\in \Ll\}$ is a family of generators of $\Mm^C$. 
\end{proposition}
\begin{proof}
For any $S\in \Ll$ we have that $E(S)$ is finite dimensional and then $E(S)^*\in \Mm^C$. It is enough to prove that the $E(S)^*$ generate the finite right comodules. If $M$ is such a comodule, then $M^*$ is a finite dimensional left comodule, thus there is a monomorphism $0\rightarrow M^*\stackrel{u}{\rightarrow}\bigoplus\limits_{\alpha\in F}E(S_\alpha)$ with $F$ a finite set and $S_\alpha$ simple left comodules. Taking duals, we obtain an epimorphism $\bigoplus\limits_{\alpha\in F}E(S_\alpha)^*\stackrel{u^*}{\rightarrow}M^{**}\cong M\rightarrow 0$ in $\Mm^C$.
\end{proof}

\begin{proposition}\label{s4}
Let $E(T)$ be an infinite dimensional injective indecomposable right comodule. Suppose that there is an epimorphism $E\stackrel{\pi}{\rightarrow} E(T)\rightarrow 0$, such that $E=\bigoplus\limits_{\alpha \in A}E_\lambda$ and $E_\lambda$ are finite dimensional injective right comodules. Then there is an epimorphism from a direct sum of finite dimensional injective right comodules to $E(T)$ with kernel containing no non-zero injective comodules. 
\end{proposition}
\begin{proof}
Denote $H={\rm Ker}\,\pi$ and consider the set $\Nn=\{Q\subset H\mid Q \,{\rm is}\,{\rm an}\,{\rm injective}\,{\rm comodule}\}$. We see that $\Nn\neq \emptyset$ as $0\in \Nn$ and that $\Nn$ is an inductive ordered set. To see this consider a chain $(X_{i})_{i\in L}$ of elements of $\Nn$ and $X=\bigcup\limits_{i\in L}X_i$ which is a subcomodule of $H$. Let $s(X)=\bigoplus\limits_{\lambda\in\Lambda}S_\lambda$ be a decomposition into simple subcomodules of the socle of $X$. Then $s(X)$ is essential in $X$ and for every $\lambda\in\Lambda$ there is an $i=i(\lambda)\in L$ such $S_\lambda\subset X_{i(\lambda)}$. As $X_i$ is injective, there is an injective envelope $H_\lambda$ of $S_\lambda$ that is contained in $X_i$. \\
First we prove that the sum $\sum\limits_{\lambda\in\Lambda}H_\lambda$ is direct. To see this it is enough to prove that $H_{\lambda_0}\cap(\sum\limits_{\lambda\in F}H_\lambda)=0$, for every finite subset $F\subseteq \Lambda$ and $\lambda_0\in\Lambda\setminus F$. We prove this by induction on the cardinal of $F$. If $F=\{\lambda\}$ then $H_{\lambda_0}\cap H_\lambda=0$ because otherwise we would have $S_{\lambda}=S_{\lambda_0}$, a contradiction. 
If the statement is proved for all sets with at most $n$ elements and $F$ is a set with $n+1$ elements then the sum $\sum\limits_{\lambda\in F}H_\lambda$ is direct, because $F=(F\setminus\{\lambda'\})\cup\{\lambda'\}$ for every $\lambda'\in F$ and we apply the induction hypothesis. If $H_{\lambda_0}\cap(\sum\limits_{\lambda\in F}H_\lambda)\neq 0$ we get that $S_{\lambda_0}\subseteq \sum\limits_{\lambda\in F}H_\lambda$, because $S_{\lambda_0}$ is essential in $H_{\lambda_0}$. But as the sum $\sum\limits_{\lambda\in F}H_\lambda$ is direct we have that $s(\sum\limits_{\lambda\in F}H_\lambda)=s(\bigoplus\limits_{\lambda\in F}H_\lambda)=\bigoplus\limits_{\lambda\in F}s(H_\lambda)=\bigoplus\limits_{\lambda\in F}S_\lambda$ so $S_{\lambda_0}\subset s(\sum\limits_{\lambda\in F}H_\lambda)=\bigoplus\limits_{\lambda\in F}S_\lambda$ which is a contradiction with $\lambda_0\notin F$. \\
Now notice that $X=\bigoplus\limits_{\lambda\in\Lambda}H_\lambda$. Since $\bigoplus\limits_{\lambda\in\Lambda}H_\lambda$ is injective, it is a direct summand of $X$.
Write $X=(\bigoplus\limits_{\lambda\in\Lambda}H_\lambda)\oplus H'$ and suppose $H'\neq 0$. Take $S'\subseteq H'$ a simple subcomodule of $H'$. Then $S'\subseteq s(X)=\bigoplus\limits_{\lambda\in\Lambda}S_\lambda\subseteq\bigoplus\limits_{\lambda\in\Lambda}H_\lambda$ which is a contradiction. We conclude that $X$ is injective, thus $X\in\Nn$.\\
By Zorn's Lemma we can then take $M$ a maximal element of $\Nn$. As $M$ is an injective comodule, it is a direct summand of $H$ and take $M\oplus H'=H$. It is obvious that $H$ is essential in $E=\bigoplus\limits_{\alpha\in A}E_\alpha$, because otherwise taking $E(H)$ an injective envelope of $H$ contained in $E$, we would have $E(H)\oplus Q=E$ so $E(T)\cong \frac{E(H)\oplus Q}{H}\cong \frac{E(H)}{H}\oplus Q$ which is a contradiction as $Q$ is a direct sum of finite dimensional comodules and $E(T)$ is indecomposable infinite dimensional. Take $E'$ an injective envelope of $H'$ contained in $E$. If $M\oplus E'\subsetneq E$ then there is a simple comodule $S$ contained in $E$ and such that $S\cap (M\oplus E')=0$, because $M\oplus E'$ is a direct summand of $E$ as it is injective. Then $S\cap H=0$ ($H\subseteq M\oplus E'$), which contradicts the fact that $H\subseteq E$ is an essential extension. Consequently, $M\oplus E'=E$ and then 
$$E(T)\cong\frac{E}{H}=\frac{M\oplus E'}{M\oplus H'}=\frac{E'}{H'}$$
where $E'$ is a direct sum of finite dimensional injective indecomposable modules and $H'$ does not contain non-zero injective modules because of the maximality of $M$.
\end{proof} 

Recall from \cite{T} that a left $C$-comodule $M$ is called quasi-finite if and only if ${\rm Hom}^C(S,M)$ is finite dimensional for every $S\in\Ll$, equivalently, if $s(M)=\bigoplus\limits_{l\in L}M_l$ is a decomposition of $M$ into simple left comodules then the set $\{l\in L:M_l\cong S\}$ is finite for every $S\in \Ll$.

\begin{lemma}\label{s5}
Let $(X_i)_{i\in L}$ be a family of nonzero (right) $C$-comodules such that $\Sigma=\bigoplus\limits_{i\in L}X_i$ is a quasifinite module. Then $\bigoplus\limits_{i\in L}X_i=\prod^{C}\limits_{i\in L}X_i$, where $\prod^C\limits_{i\in L}$ is the product in the category of comodules.
\end{lemma}
\begin{proof}
We have that $\prod^C\limits_{i\in L}X_i=Rat(\prod\limits_{i\in L}X_i)$, where $\prod\limits_{i\in L}$ is the product of modules. Suppose that $x=(x_i)_{i\in L}\in P=Rat(\prod\limits_{i\in L}X_i)$ and the set $L^\prime=\{i\in L\mid x_i\neq 0\}$ is infinite. Then $C^*\cdot x$ is a finite dimensional rational module, so it has a finite composition series. For each $i\in L$ the canonical projection $p_i:C^*\cdot x\rightarrow C^*\cdot x_i$ is an epimorphism, thus $C^*\cdot x_i$ is a rational module. 
For every $i\in L^\prime$ consider $S_i$ a simple subcomodule of $C^*\cdot x_i$. As $\Sigma$ is quasi-finite and $L^\prime$ is infinite, we have that the set $\Rr'=\{T\in \Rr\mid \exists\,i\in L^\prime\,{\rm such}\,{\rm that}\,S_i\cong T\}$ is infinite ($\Rr$ is the chosen set of representatives for the types of isomorphisms of simple right $C$-comodules). For each $T\in \Rr'$ choose $k\in L^\prime$ such that $T\cong S_k$. Denote by $\Lambda$ the set of these $k$'s. As for every $T\in \Rr'$ there is $k$, a monomorphism $T\hookrightarrow C^*\cdot x_k$ and an epimorphism $C^*\cdot x\stackrel{p_k}{\longrightarrow} C^*\cdot x_k$, it follows then that every composition series of $C^*\cdot x$ contains a simple factor isomorphic to $T$. As $C^*\cdot x$ is finite dimensional it follows that the set $\Ff$ of simple left $C^*$-modules appearing as factors in any composition series is finite. But $\Rr'\subseteq \Ff$ which is a contradiction to the fact that $\Rr'$ is infinite. Thus $x\in\bigoplus\limits_{i\in L}X_i$, and then $Rat(\prod\limits_{i\in L}X_i)\subseteq\bigoplus\limits_{i\in L}X_i$ and the proof is finished as the converse inclusion is obviously true.
\end{proof}

\begin{theorem}\label{s6}
Let $C$ be a coalgebra such that $C\cong Rat(C^*_{C^*})$ as right $C^*$-modules (left $C$-comodules). Then $C$ is left semiperfect and $C\cong Rat({}_{C^*}C^*)$ as left $C^*$-modules. 
\end{theorem}
\begin{proof}
Let $T$ be a simple right $C$-comodule such that $E(T)$ is infinite dimensional. Then $Rat(E(T)^*)=0$ by proposition \ref{s1}. We have that $C$ is right co-Frobenius, thus it is also right semiperfect. By Proposition \ref{s3} there is an exact sequence of right comodules $$0\rightarrow H\hookrightarrow \bigoplus\limits_{\alpha\in A}E(S_\alpha)^*\stackrel{u}{\rightarrow}E(T)\rightarrow 0$$
with $E(S_\alpha)$ finite dimensional injective left comodules. Let $E=\bigoplus\limits_{\alpha\in A}E(S_\alpha)^*$.  As $E(S_\alpha)$ are finite dimensional injective left comodules, they are injective also as right $C^*$-modules (by Lemma \ref{i1}) and as $C\cong Rat(C^*_{C^*})$ it follows that every injective indecomposable comodule is a direct summand of $C^*$, thus it is projective. 
By Lemma \ref{i2} we obtain that every $E_\alpha=E(S_\alpha)^*$ is also injective, and also finite dimensional indecomposable and then by Proposition \ref{s4} we may assume that $H$ does not contain nonzero injective comodules. Take $n\in A$ and set $E'=\bigoplus\limits_{\alpha\in A\setminus \{n\}}E_\alpha$. Then $H+E'=E$. To see this first notice that $H+E'$ has finite codimension in $E$. 
There is an epimorphism of right comodules (thus of left $C^*$-modules) 
$$E(T)\cong\frac{E}{H}\rightarrow\frac{E}{H+E'}\rightarrow 0$$ 
and by taking duals we get a monomorphism of right $C^*$-modules 
$$0\rightarrow \left(\frac{E}{H+E'}\right)^*\rightarrow\left(\frac{E}{H}\right)^*\cong E(T)^*$$
But the dual of the finite dimensional right comodule $E/(H+E')$ is a rational right $C^*$-module, implying that $Rat(E(T)^*)\neq 0$ if $H+E'\neq E$, a contradiction. \\
By the isomorphisms
$$E_n\cong \frac{E}{E'}=\frac{H+E'}{E'}\cong \frac{H}{H\cap E'}$$
we conclude that there is an epimorphism from $H$ onto $E_n$. This morphism must split, as $E_n$ is also a projective right comodule (again by Lemma \ref{i2}). This shows that $H$ contains an injective subcomodule isomorphic to $E_n$, which contradicts the supposition that $H$ does not contain nonzero injective subcomodules. Hence $C$ is semiperfect.\\
Now we have $\bigoplus\limits_{i\in I}E(S_i)\cong Rat(C^*_{C^*})\cong Rat(\prod\limits_{j\in J_0}E(T_j)^*)$. But $E(T_j)^*$ are indecomposable and also injective finite dimensional left comodules; write $L_j=E(T_j)^*$. Then $L_j$ is the injective envelope of its socle, so $L_j\cong L_{j'}$ if and only if $s(L_j)\cong s(L_{j'})$ (and equivalently, $E(T_j)=L_j^*\cong L_{j'}^*=E(T_{j'})$). This shows that for any $S\in \Ll$, there are only finitely many $j\in J$ with the property $s(L_j)\cong S$, because only finitely many $E(T_j)$'s can be isomorphic to the same injective indecomposable. This shows that the comodule $\bigoplus\limits_{j\in J}L_j$ is quasifinite, and then by Lemma \ref{s5} we have that $\bigoplus\limits_{i\in I}E(S_i)\cong Rat(\prod\limits_{j\in J}E(T_j)^*)\cong\bigoplus\limits_{j\in J}E(T_j)^*$. By Krull-Remak-Schmidt-Azumaya's Theorem, there is a bijection $\varphi:I\rightarrow J$ such that $E(S_i)\cong E(T_{\varphi(i)})^*$ for every $i\in I$, equivalently, $E(S_i)^*\cong E(T_{\varphi(i)})$ for all $i$. We then obtain that the comodule $\bigoplus\limits_{j\in J}E(T_j)\cong\bigoplus\limits_{i\in I}E(S_i)^*$ is quasifinite, and again by Lemma \ref{s5} we have that $\bigoplus\limits_{j\in J}E(T_j)\cong Rat(\prod\limits_{i\in I}E(S_i)^*)$, finally showing that $C\cong Rat({}_{C^*}C^*)$ as left $C^*$-modules.
\end{proof}

\begin{remark}
Let $\Cc$ be a category. Then for every two objects $X,Y$, we consider the "Yoneda" bijection of sets $\Lambda:\underline{\underline{\mathbf{ Nat}}}({\rm Hom}(-,X),{\rm Hom}(-,Y))\longrightarrow{\rm Hom}(X,Y)$ defined by $\Lambda(\varphi)=\varphi_X(1_X)$ with inverse $\Lambda^{-1}(\theta)=(f\mapsto \theta\circ f)$. 
Moreover, if $\varphi$ is a natural equivalence with inverse $\varphi'$, then $\theta=\varphi_X(1_X):X\rightarrow Y$ is an isomorphism with inverse $\theta'=\varphi'_Y(1_Y):Y\rightarrow X$.
\end{remark}


\begin{theorem}\label{s8}
Let $C$ be a coalgebra. Then the following assertions are equivalent:
\begin{enumerate}
\item[(i)] $C$ is a co-Frobenius coalgebra (left and right co-Frobenius).
\item[(ii)] $Rat({}_{C^*}C^*)\cong C$ in $\Mm^C$ (or as left $C^*$-modules). 
\item[(iii)] The functors ${\rm Hom}_K(-,K)$ and ${\rm Hom}_{C^*}(-,C^*)$ from $\Mm^C\subset {}_{C^*}\Mm$ to $\Mm_{C^*}$ are naturally equivalent.
\item[(iv)] The left hand side versions of (ii) and (iii).
\item[(v)] There is a $K$-bilinear form $(-,-)$ on $C\times C$ that is $C^*$-balanced (i.e. $(c\cdot h^*,d)=(c,h^*\cdot d)$ for all $c,d\in C$ and $h^*\in C^*$) and left and right non-degenerate. 
\end{enumerate}
\end{theorem}
\begin{proof}
(ii)$\Leftrightarrow$(iii) We have a natural equivalence $\rm{Hom}^C(-,C)\cong\rm{Hom}(-,K):\Mm^C\rightarrow\Mm_{C^*}$, $h\mapsto\varepsilon\circ h$, where $\varepsilon$ is the counit of $C$ and also a natural equiavalence of functors ${\rm Hom}_{{}_{C^*}\Mm}(-,C^*)\cong {\rm Hom}^C(-,C^*)\cong{\rm Hom}^{C}(-,Rat({}_{C^*}C^*))$. Thus by the previous remark the functors ${\rm Hom}(-,K)$ and ${\rm Hom}_{C^*}(-,C^*)$ from $\Mm^C$ to $\Mm_{C^*}$ are isomorphic if and only if there is an isomorphism of left $C^*$-modules (right $C$-comodules) $C\cong Rat({}_{C^*}C^*)$. \\   
(i)$\Rightarrow$(ii) follows from Proposition \ref{i3} and (ii)$\Rightarrow$(i) from Theorem \ref{s6}.\\
(v)$\Leftrightarrow$(i) For an isomorphism $\varphi:C\rightarrow Rat({}_{C^*}C^*)$ one can define $(-,-):C\times C\rightarrow C$ by $(c,d)=\varphi(d)(c)$. This is a $C^*$-balanced form by the morphism property of $\varphi$, is left non-degenerate by the injectivity of $\varphi$ and right nondegenerate by the density of ${\rm Im}\varphi=Rat({}_{C^*}C^*)$ in $C^*$ in the finite topology on $C^*$, because $Rat({}_{C^*}C^*)=\bigoplus\limits_{i\in I}E(S_i)^*$ (as shown above in the proof of Theorem \ref{s6}) which is dense in $\prod\limits_{i\in I}E(S_i)^*$. For the converse consider $\varphi:C\rightarrow C^*$ and $\psi:C\rightarrow C^*$ defined by $\varphi(c)(d)=(c,d)$ and $\psi(c)(d)=(d,c)$ for all $c,d\in C$. Then an easy computation shows that $\varphi$ is an injective morphism of right $C^*$-modules and $\psi$ is an injective morphism of left comodules, thus $C$ is left and right co-Frobenius. 
\end{proof}

\begin{center}
\sc Acknowledgment
\end{center}
The author wishes to thank his Ph.D. advisor professor Constantin N\u ast\u asescu for very useful remarks on the subject as well as for his continuous support throughout the past years. He also wishes to thank the referee for very useful and helpful remarks that greately improved the presentation of this paper.

\vspace*{3mm} 
\begin{flushright}
\begin{minipage}{148mm}\sc\footnotesize

Miodrag Cristian Iovanov\\
University of Bucharest, Faculty of Mathematics, Str.
Academiei 14,
RO-70109, Bucharest, Romania\\
{\it E--mail address}: {\tt
yovanov@gmail.com}\vspace*{3mm}

\end{minipage}
\end{flushright}
\end{document}